\documentclass[11pt]{article}
\usepackage[russian]{babel}
\usepackage[cp1251]{inputenc}
\usepackage{amsmath,latexsym}
\usepackage[psamsfonts]{amssymb}
\usepackage[dvips]{graphicx}
\usepackage[mathcal]{euscript}
\begin{document}
\begin{center}
\textbf{ Построение фундаментальной системы решений \\ для вырождающегося уравнения с дробной производной Джрбашяна-Нерсесяна }\\
\textbf{Б.Ю.Иргашев}\\
Наманганский инженерно-строительный институт, Узбекистан\\
Институт Математики им.В.И.Романовского АН РУз.\\
bahromirgasev@gmail.com
\end{center}
УДК.517.926.4

\textbf{Аннотация.} В статье построено общее решение одного вырождающегося уравнения с дробной производной Джрбашяна-Нерсесяна. Частные решения представлены через функцию Килбаса-Сайго.

 \textbf{Kлючевые слова.} Производная дробного порядка , вырождение, ряд, функция Килбаса-Сайго, решение.

В последнее время специалистами интенсивно изучаются уравнения c участием производных дробного порядка с переменными коэффициентами. К числу таких уравнений относятся вырождающиеся уравнения.  В работе  [1] изучалось уравнение
\[D_{0x}^\alpha {t^\beta }u\left( t \right) = \lambda u\left( x \right),\,\,0 < x < b,\]	
где $0 < \alpha  < 1,\,\,\lambda  - $ спектральный параметр, $\beta  = const \ge 0.$
В работе [2] были найдены решения в замкнутой форме уравнений дробного порядка
\[\left( {D_{0 + }^\alpha y} \right)\left( x \right) = a{x^\beta }y\left( x \right) + f\left( x \right)\left( {0 < x < d \le \infty ,\alpha  > 0,\beta  \in R,a \ne 0} \right),\]
\[\left( {D_ - ^\alpha y} \right)\left( x \right) = a{x^\beta }y\left( x \right) + f\left( x \right)\left( {0 \le d < x < \infty ,\alpha  > 0,\beta  \in R,a \ne 0} \right),\]
с дробными производными Римана-Лиувилля на полуоси $\left( {0,\infty } \right)$ [3].\\
 К таким уравнениям приводят прикладные задачи [4]. Пример такого уравнения дает уравнение теории полярографии [5]
\[\left( {D_{0 + }^{{1 \mathord{\left/
 {\vphantom {1 2}} \right.
 \kern-\nulldelimiterspace} 2}}y} \right)\left( x \right) = a{x^\beta }y\left( x \right) + {x^{ - {1 \mathord{\left/
 {\vphantom {1 2}} \right.
 \kern-\nulldelimiterspace} 2}}},\left( {0 < x, - {1 \mathord{\left/
 {\vphantom {1 {2 < \beta  \le }}} \right.
 \kern-\nulldelimiterspace} {2 < \beta  \le }}0} \right),\]
 возникающее при $a =  - 1$  в задачах диффузии [5].

Рассмотрим следующее уравнение
\[D_{0y}^{\left\{ {{\gamma _0},{\gamma _1},...,{\gamma _{m - 1}},{\gamma _m}} \right\}}u\left( y \right) = \lambda {y^s}u(y),y > 0,\lambda  \in C,s \ge 0,\eqno(1)\]
где $D_{0y}^{\left\{ {{\gamma _0},{\gamma _1},...,{\gamma _{m - 1}},{\gamma _m}} \right\}}$ - оператор дробного дифференцирования Джрбашяна–Нерсесяна порядка $\alpha  = \sum\limits_{k = 0}^m {{\gamma _k}}  - 1 > 0,$  ассоциированный с последовательностью $ \left\{ {{\gamma _k}} \right\}_0^m = \left\{ {{\gamma _0},{\gamma _1},...,{\gamma _{m - 1}},{\gamma _m}} \right\}$ , ${\gamma _k} \in \left( {0,1} \right],k = 0,1,...,m,$ определяется соотношением [6]
\[D_{0y}^{\left\{ {{\gamma _0},{\gamma _1},...,{\gamma _{m - 1}},{\gamma _m}} \right\}} = D_{0y}^{{\gamma _m} - 1}D_{0y}^{{\gamma _{m - 1}}}...D_{0y}^{{\gamma _1}}D_{0y}^{{\gamma _0}},\eqno (2)\]
здесь $D_{0y}^\gamma $ - оператор дробного интегро-дифференцирования в смысле Римана-Лиувилля порядка $\gamma$ с началом в точке $ y=0 $ определяемый следующим образом [1, с. 9]
\[D_{0y}^\gamma g\left( y \right) = \left\{ \begin{array}{l}
\frac{1}{{\Gamma \left( { - \gamma } \right)}}\int\limits_0^y {\frac{{g\left( t \right)dt}}{{{{\left| {y - t} \right|}^{1 + \gamma }}}}} ,\gamma  < 0,\\
g\left( y \right),\gamma  = 0,\\
{\left( {\frac{d}{{dy}}} \right)^p}D_{0y}^{\gamma  - p}g\left( y \right),p - 1 < \gamma  \le p,p \in N.
\end{array} \right.\]
Заметим, что если в (1) в качестве последовательности $\left\{ {{\gamma _k}} \right\}_0^m$ взять последовательность $\left\{ {{\gamma _k}} \right\}_0^m = \left\{ {\alpha  - m + 1,\underbrace {1,...,1}_m} \right\},$ то мы получим производную Римана–Лиувилля:
\[D_{0y}^{\left\{ {\alpha  - m + 1,1,...,1} \right\}} = D_{0y}^\alpha ,m - 1 < \alpha  \le m.\]
Последовательности $\left\{ {{\gamma _k}} \right\}_0^m = \left\{ {\underbrace {1,...,1,}_m\alpha  - m + 1,} \right\},$ соответствует производная Капуто:
\[D_{0y}^{\left\{ {1,...,1,\alpha  - m + 1} \right\}} = {}_CD_{0y}^\alpha ,m - 1 < \alpha  \le m.\]
В работе [6] рассматривалась задача Коши для уравнения вида
\[\sum\limits_{k = 0}^n {{a_k}D_{0y}^{\left\{ {{\gamma _0},...,{\gamma _k}} \right\}}u\left( y \right)}  = f\left( y \right),\eqno(3)\]
с переменными коэффициентами. Исследуемая задача эквивалентно сведена к интегральному уравнению Вольтерра второго рода. Доказана теорема существования и единственности решения. В работе [7] в терминах функции Райта строится явное представление решения задачи Коши для уравнения (3). В работе [8] для линейного обыкновенного дифференциального уравнения дробного порядка вида (3) с производными Римана-Лиувилля была сформулирована и решена начальная задача. Краевые и начальные задачи для вырождающихся уравнений с дробным производным Хилфера исследовались в работах [9-12], а с дробными производными Римана-Лиувилля и Капуто в работах [2],[13-14].

В данной работе в терминах функции Килбаса-Сайго строится явное представление фундаментальной системы решений уравнения (1).

Приступим к построению решения уравнения (1). Решение будем искать в виде
\[u\left( y \right) = \sum\limits_{n = 0}^\infty  {{c_n}{y^{an + b}}},\eqno(4)\]
где ${c_n},a>0,b$ пока неизвестные вещественные числа.\\
Сделаем предварительные вычисления, имеем:
\[D_{0y}^\gamma {y^\delta } = \frac{d}{{dy}}D_{0y}^{{\gamma _0} - 1}{y^\delta } = \frac{d}{{dy}}\frac{1}{{\Gamma \left( {1 - {\gamma _0}} \right)}}\int\limits_0^y {\frac{{{t^\delta }dt}}{{{{\left( {y - t} \right)}^{{\gamma _0}}}}}}  = \frac{{\left( {\delta  + 1 - {\gamma _0}} \right)\Gamma \left( {\delta  + 1} \right){y^{\delta  - {\gamma _0}}}}}{{\Gamma \left( {\delta  - {\gamma _0} + 2} \right)}},\delta  >  - 1.\]
Далее
\[D_{0y}^{{\gamma _1}}D_{0y}^{{\gamma _0}}{y^\delta } = \frac{{\left( {\delta  + 1 - {\gamma _0}} \right)\Gamma \left( {\delta  + 1} \right)}}{{\Gamma \left( {\delta  - {\gamma _0} + 2} \right)}}D_{0y}^{{\gamma _1}}{y^{\delta  - {\gamma _0}}} = \]
\[ = \frac{{\left( {\delta  - {\gamma _0} + 1} \right)\left( {\delta  - {\gamma _0} - {\gamma _1} + 1} \right)\Gamma \left( {\delta  + 1} \right)\Gamma \left( {\delta  - {\gamma _0} + 1} \right)}}{{\Gamma \left( {\delta  - {\gamma _0} + 2} \right)\Gamma \left( {\delta  - {\gamma _0} - {\gamma _1} + 2} \right)}}{y^{\delta  - {\gamma _0} - {\gamma _1}}},\delta  - {\gamma _0} >  - 1.\]
Продолжая этот процесс получим
\[D_{0y}^{{\gamma _{m - 1}}}...D_{0y}^{{\gamma _0}}{y^\delta } = \frac{{\Gamma \left( {\delta  + 1} \right)\prod\limits_{k = 0}^{m - 1} {\left( {\delta  - {\alpha _k}} \right)} \prod\limits_{k = 0}^{m - 2} {\Gamma \left( {\delta  - {\alpha _k}} \right)} }}{{\prod\limits_{k = 0}^{m - 1} {\Gamma \left( {\delta  - {\alpha _k} + 1} \right)} }}{y^{\delta  - {\alpha _{m - 1}} - 1}},\delta  - {\alpha _{m-2}} > 0,\]
где
\[{\alpha _k} = \sum\limits_{j = 0}^k {{\gamma _j}}  - 1,{\alpha _m} = \alpha .\]
Окончательно имеем формулу
\[D_{0y}^{\left\{ {{\gamma _0},{\gamma _1},...,{\gamma _{m - 1}},{\gamma _m}} \right\}}{y^\delta } = \]
\[ = \frac{{\Gamma \left( {\delta  + 1} \right)\prod\limits_{k = 0}^{m - 1} {\left( {\delta  - {\alpha _k}} \right)} \prod\limits_{k = 0}^{m - 1} {\Gamma \left( {\delta  - {\alpha _k}} \right)} {y^{\delta  - \alpha }}}}{{\prod\limits_{k = 0}^m {\Gamma \left( {\delta  - {\alpha _k} + 1} \right)} }},\delta  - {\alpha _{m - 1}} > 0.\eqno(5)\]
Теперь подставим (4) в (1), затем используя формулу (5) получим формальное равенство
\[\sum\limits_{n = 0}^\infty  {{c_n}} \frac{{\Gamma \left( {an + b + 1} \right)\prod\limits_{k = 0}^{m - 1} {\left( {an + b - {\alpha _k}} \right)} \prod\limits_{k = 0}^{m - 1} {\Gamma \left( {an + b - {\alpha _k}} \right)} }}{{\Gamma \left( {an + b + 1 - {\alpha}} \right)\prod\limits_{k = 0}^{m - 1} {\Gamma \left( {an + b + 1 - {\alpha _k}} \right)} }}{y^{an + b - \alpha }}=\]
\[ = \lambda \sum\limits_{n = 0}^\infty  {{c_n}{y^{an + b + s}}} ,\]
Пусть
\[a = \alpha  + s,\]
\[b = {\alpha _k},k = 0,1,...,m - 1,\]
тогда используя равенство:
\[\frac{{\prod\limits_{k = 0}^{m - 1} {\left( {an + b - {\alpha _k}} \right)} \prod\limits_{k = 0}^{m - 1} {\Gamma \left( {an + b - {\alpha _k}} \right)} }}{{\prod\limits_{k = 0}^{m - 1} {\Gamma \left( {an + b + 1 - {\alpha _k}} \right)} }} = 1,\]
получим
\[\sum\limits_{n = 0}^\infty  {{c_n}} \frac{{\Gamma \left( {an + b + 1} \right)}}{{\Gamma \left( {an + b + 1 - {\alpha}} \right)}}{y^{a\left( {n - 1} \right)}} = \lambda \sum\limits_{n = 0}^\infty  {{c_n}{y^{an}}}.\]
Найдем неизвестные коэффициенты $c_{n}$
\[{c_n} = \lambda {c_{n - 1}}\frac{{\Gamma \left( {an + b + 1 - \alpha } \right)}}{{\Gamma \left( {an + b + 1} \right)}} = {\lambda ^n}{c_0}\frac{{\prod\limits_{j = 0}^{n - 1} {\Gamma \left( {aj + a + b + 1 - \alpha } \right)} }}{{\prod\limits_{j = 0}^{n - 1} {\Gamma \left( {aj + a + b + 1} \right)} }}.\]
Заметим,что
\[\begin{array}{l}
ja + a + b - \alpha  + 1 \ge a + b - \alpha  + 1 \ge \\
 \ge \alpha  + s + {\gamma _0} - 1 - \alpha  + 1 = s + {\gamma _0} \ge {\gamma _0} > 0.
\end{array}\]
Итак получили следующее семейство линейно независимых решений уравнения (1)
\[{u_k}\left( y \right) = {y^{{\alpha _k}}}\sum\limits_{n = 0}^\infty  {{c_n}{{\left( {\lambda {y^{\alpha  + s}}} \right)}^n}} ,k = 0,1,...,m - 1,\eqno(6)\]
где
\[{c_0} = 1,{c_n} = \frac{{\prod\limits_{j = 0}^{n - 1} {\Gamma \left( {\left( {\alpha  + s} \right)j + s + b + 1} \right)} }}{{\prod\limits_{j = 0}^{n - 1} {\Gamma \left( {\left( {\alpha  + s} \right)j + s + b + 1 + \alpha } \right)} }},n = 1,2,....\eqno (7)\]
Покажем абсолютную сходимость ряда (6). Применим признак Даламбера, имеем
\[{y^{\alpha  + s}}\mathop {\lim }\limits_{n \to  + \infty } \left\{ {\frac{{\prod\limits_{j = 0}^n {\Gamma \left( {\left( {\alpha  + s} \right)j + s + b + 1} \right)} }}{{\prod\limits_{j = 0}^n {\Gamma \left( {\left( {\alpha  + s} \right)j + s + b + 1 + \alpha } \right)} }}:\frac{{\prod\limits_{j = 0}^{n - 1} {\Gamma \left( {\left( {\alpha  + s} \right)j + s + b + 1} \right)} }}{{\prod\limits_{j = 0}^{n - 1} {\Gamma \left( {\left( {\alpha  + s} \right)j + s + b + 1 + \alpha } \right)} }}} \right\} = \]
\[ = {y^{\alpha  + s}}\mathop {\lim }\limits_{n \to  + \infty } \frac{{\Gamma \left( {\left( {\alpha  + s} \right)n + s + b + 1} \right)}}{{\Gamma \left( {\left( {\alpha  + s} \right)n + s + b + 1 + \alpha } \right)}} = {y^{\alpha  + s}}\mathop {\lim }\limits_{n \to  + \infty } {\left\{ {\left( {\alpha  + s} \right)n} \right\}^{ - \alpha }} = 0,\]
т.к. [15]
\[\frac{{\Gamma \left( {z + \alpha } \right)}}{{\Gamma \left( {z + \beta } \right)}} = O\left( {{z^{\alpha  - \beta }}} \right),z \to  + \infty .\]
Представление (7) также можно записать в виде
\[\frac{{\prod\limits_{j = 0}^{n - 1} {\Gamma \left( {\left( {\alpha  + s} \right)j + s + b + 1} \right)} }}{{\prod\limits_{j = 0}^{n - 1} {\Gamma \left( {\left( {\alpha  + s} \right)j + s + b + 1 + \alpha } \right)} }} = \frac{{\prod\limits_{j = 0}^{n - 1} {\Gamma \left( {\alpha \left( {\frac{{\alpha  + s}}{\alpha }j + \frac{{s + {\alpha _k}}}{\alpha }} \right) + 1} \right)} }}{{\prod\limits_{j = 0}^{n - 1} {\Gamma \left( {\alpha \left( {\frac{{\alpha  + s}}{\alpha }j + \frac{{s + {\alpha _k}}}{\alpha } + 1} \right) + 1} \right)} }},\]
тогда семейство линейно независимых решений запишется так
\[{u_k}\left( y \right) = {y^{{\alpha _k}}}{E_{\alpha ,\frac{{\alpha  + s}}{\alpha },\frac{{{\alpha _k} + s}}{\alpha }}}\left( {\lambda {y^{\alpha  + s}}} \right),k = 0,1,...,m - 1,\eqno (8)\]
где
\[{E_{\alpha ,m,l}}\left( z \right) = \sum\limits_{i = 0}^\infty  {{c_i}{z^i}} ,{c_0} = 1,{c_i} = \prod\limits_{j = 0}^{i - 1} {\frac{{\Gamma \left( {\alpha \left( {jm + l} \right) + 1} \right)}}{{\Gamma \left( {\alpha \left( {jm + l + 1} \right) + 1} \right)}}} ,i \ge 1\]
- функция Килбаса-Сайго (см.[2]).
Итак общее решение уравнения (1) имеет вид
\[u\left( y \right) = \sum\limits_{k = 0}^{m - 1} {{d_k}{y^{{\alpha _k}}}{E_{\alpha ,\frac{{\alpha  + s}}{\alpha },\frac{{{\alpha _k} + s}}{\alpha }}}\left( {\lambda {y^{\alpha  + s}}} \right)} ,{d_k} = const,k = 0,1,...,m - 1.\eqno(9)\]
Рассмотрим некоторые частные случаи.\\
1. Пусть в уравнении (1) $s=0$, тогда из формулы [2]
\[{E_{\alpha ,1,l}}\left( z \right) = \Gamma \left( {\alpha l + 1} \right){E_{\alpha ,\alpha l + 1}}\left( {\lambda {y^\alpha }} \right),\]
следует
\[{u_k}\left( y \right) = {y^{{\alpha _k}}}{E_{\alpha ,1,\frac{{{\alpha _k}}}{\alpha }}}\left( {\lambda {y^\alpha }} \right) = \Gamma \left( {{\alpha _k} + 1} \right){E_{\alpha ,{\alpha _k} + 1}}\left( {\lambda {y^\alpha }} \right),\]
где
\[{E_{\alpha ,\beta }}\left( z \right) = \sum\limits_{n = 0}^\infty  {\frac{{{z^n}}}{{\Gamma \left( {\alpha n + \beta } \right)}}} ,\left( {\alpha ,\beta  > 0} \right)\]
- функция Миттаг-Леффлера.
Что с точностью до множителя совпадает с результатами из работы [6, представление (3.15)].

2. Пусть в уравнении (1) имеем оператор Римана-Лиувилля , т.е.
\[D_{0y}^{\left\{ {\alpha  - m + 1,\underbrace {1,...,1}_m} \right\}}u\left( y \right) = D_{0y}^\alpha u\left( y \right) = \lambda {y^s}u,\,m - 1 < \alpha  \le m,m \in N,\]
тогда учитывая, что
\[{\alpha _k} = {\gamma _0} + ... + {\gamma _k} - 1 = \alpha  + k - m,k = 0,1,...,m - 1,\]
из (8) имеем
\[{u_k}\left( y \right) ={y^{\alpha  + k - m}}{E_{\alpha ,\frac{{\alpha  + s}}{\alpha },\frac{{\alpha  + k - m + s}}{\alpha }}}\left( {\lambda {y^{\alpha  + s}}} \right) = \]
\[ = {y^{\alpha  - \left( {m - k} \right)}}{E_{\alpha ,\frac{{\alpha  + s}}{\alpha },\frac{{\alpha  + s - \left( {m - k} \right)}}{\alpha }}}\left( {\lambda {y^{\alpha  + s}}} \right) = \left( {j = m - k} \right)\]
\[ = {y^{\alpha  - j}}{E_{\alpha ,\frac{{\alpha  + s}}{\alpha },\frac{{\alpha  + s - j}}{\alpha }}}\left( {\lambda {y^{\alpha  + s}}} \right),j = 1,2,...,m.\eqno(10)\]
Представление (10) совпадает с результатами из работы [2, формулы (19),(21)].

3. Пусть в уравнении (1) имеем оператор Капуто , т.е.
\[D_{0y}^{\left\{ {\underbrace {1,...,1,}_m\alpha  - m + 1} \right\}}u\left( y \right) = \lambda {y^s}u\left( y \right),\]
далее имеем
\[D_{0y}^{\left\{ {\underbrace {1,...,1,}_m\alpha  - m + 1} \right\}} = D_{0y}^{\alpha  - m}{\left( {\frac{d}{{dy}}} \right)^m},\]
\[{\alpha _0} = {\gamma _0} - 1 = 0,{\alpha _1} = {\gamma _0} + {\gamma _1} - 1 = 1,\]
\[{\alpha _k} = k,k = 0,1,...,m - 1,\]
\[{\alpha _m} = m + \alpha  - m + 1 - 1 = \alpha ,\]
отсюда
\[{u_k}\left( y \right) = {y^{{\alpha _k}}}{E_{\alpha ,\frac{{\alpha  + s}}{\alpha },\frac{{{\alpha _k} + s}}{\alpha }}}\left( {\lambda {y^{\alpha  + s}}} \right) = \]
\[ = {y^k}{E_{\alpha ,\frac{{\alpha  + s}}{\alpha },\frac{{k + s}}{\alpha }}}\left( {\lambda {y^{\alpha  + s}}} \right),k = 0,1,...,m - 1.\]
 Это овпадает с результатами из работы [14, формула (4.1.82)].

4. Пусть в уравнение (1) имеем дробный оператор Хилфера [11]:
\[I_{0y}^{\mu \left( {m - \alpha } \right)}{\left( {\frac{d}{{dt}}} \right)^m}I_{0y}^{\left( {1 - \mu } \right)\left( {m - \alpha } \right)}u\left( y \right) = \lambda {y^s}u\left( y \right),0 \le \mu  \le 1,m - 1 < \alpha  \le m,y > 0.\]
Оператор Хилфера запишем в виде оператора Джрбашяна-Нерсесяна
\[I_{0y}^{\mu \left( {m - \alpha } \right)}{\left( {\frac{d}{{dt}}} \right)^m}I_{0y}^{\left( {1 - \mu } \right)\left( {m - \alpha } \right)} = D_{0y}^{1 - \mu \left( {m - \alpha } \right) - 1}{\left( {\frac{d}{{dt}}} \right)^{m - 1}}D_{0y}^{1 - \left( {1 - \mu } \right)\left( {m - \alpha } \right)} = \]
\[ = D_{0y}^{\left\{ {1 - \left( {1 - \mu } \right)\left( {m - \alpha } \right),\underbrace {1,...,1}_{m - 1},1 - \mu \left( {m - \alpha } \right)} \right\}},\]
отсюда имеем
\[{\alpha _k} =  - \left( {1 - \mu } \right)\left( {m - \alpha } \right) + k,k = 0,1,...,m - 1,\]
\[{\alpha _m} =  - \left( {1 - \mu } \right)\left( {m - \alpha } \right) + m - 1 + 1 - \mu \left( {m - \alpha } \right) =  - \left( {m - \alpha } \right) + m = \alpha .\]
Теперь из представления (8) имеем
\[{u_k}\left( y \right) = {y^{{\alpha _k}}}{E_{\alpha ,\frac{{\alpha  + s}}{\alpha },\frac{{{\alpha _k} + s}}{\alpha }}}\left( {\lambda {y^{\alpha  + s}}} \right) = {y^{k - \left( {1 - \mu } \right)\left( {m - \alpha } \right)}}{E_{\alpha ,\frac{{\alpha  + s}}{\alpha },\frac{{s + k - \left( {1 - \mu } \right)\left( {m - \alpha } \right)}}{\alpha }}}\left( {\lambda {y^{\alpha  + s}}} \right).\]

Применим формулу (9) к получению представления решения следующей задачи Коши (см.[6]):
\[\left\{ \begin{array}{l}
D_{0y}^{\left\{ {{\gamma _0},{\gamma _1},...,{\gamma _{m - 1}},{\gamma _m}} \right\}}u\left( y \right) = \lambda {y^s}u,y > 0,\lambda  \in C,s \ge 0.\\
\mathop {\lim }\limits_{y \to 0} D_{0y}^{{\alpha _0}}u\left( y \right) = {A_0},\\
\mathop {\lim }\limits_{y \to 0} D_{0y}^{{\alpha _1}}u\left( y \right) = {A_1},\\
...\\
\mathop {\lim }\limits_{y \to 0} D_{0y}^{{\alpha _{m - 1}}}u\left( y \right)u\left( y \right) = {A_{m - 1}}.
\end{array} \right.\eqno(11)\]
здесь
\[{A_i} = const,i = 0,1,...,m - 1,\]
\[D_{0y}^{{\alpha _0}} = D_{0y}^{{\gamma _0} - 1},\]
\[D_{0y}^{{\alpha _k}} = D_{0y}^{{\gamma _k} - 1}\frac{d}{{dy}}D_{0y}^{{\alpha _{k - 1}}}.\]
Справедлива формула [6]
\[D_{0y}^{{\alpha _s}}{y^{{\alpha _k}}} = \left\{ \begin{array}{l}
0,0 \le k \le s - 1,\\
\Gamma \left( {1 + {\alpha _k}} \right),k = s,\\
\frac{{\Gamma \left( {1 + {\alpha _k}} \right)}}{{\Gamma \left( {1 + {\alpha _k} - {\alpha _s}} \right)}}{y^{{\alpha _k} - {\alpha _s}}},s < k \le m.
\end{array} \right.\]
Подставив представление (9) в начальные условия (11) , получим
\[\mathop {\lim }\limits_{y \to 0} D_{0y}^{{\alpha _0}}u\left( y \right) = {d_0}\Gamma \left( {1 + {\alpha _0}} \right) = {A_0},\]
\[\mathop {\lim }\limits_{y \to 0} D_{0y}^{{\alpha _1}}u\left( y \right) = {d_1}\Gamma \left( {1 + {\alpha _1}} \right) = {A_1},\]
\[...\]
\[\mathop {\lim }\limits_{y \to 0} D_{0y}^{{\alpha _{m - 1}}}u\left( y \right) = {d_{m - 1}}\Gamma \left( {1 + {\alpha _{m - 1}}} \right) = {A_{m - 1}}.\]
Значит решение задачи Коши будет иметь вид
\[u\left( y \right) = \sum\limits_{k = 0}^{m - 1} {\frac{{{A_k}{y^{{\alpha _k}}}}}{{\Gamma \left( {{\alpha _k} + 1} \right)}}{E_{\alpha ,\frac{{\alpha  + s}}{\alpha },\frac{{{\alpha _k} + s}}{\alpha }}}\left( {\lambda {y^{\alpha  + s}}} \right)} .\]

\begin{center}
Литература
\end{center}
1. Нахушев А.М. Дробное исчисление и его применение. - М.: Физматлит. 2003. - 272 c.\\
2. Килбас А. А., Сайго М.  Решение в замкнутой форме одного класса линейных дифференциальных уравнений дробного порядка. Дифференц. уравнения, 33 (2), 1997. c. 195 - 204.\\
3. Самко С. Г., Килбас А. А., Маричев О. И. Интегралы и производные дробного порядка и некоторые их приложения. Минск. Наука и техника. 1987. - 688 с.\\
4. Oldham К. В., Spanier J. The fractional calculus. New York; London. 1974.\\
5. Wiener K.  Wiss. Z. Univ. Halle Math. Natur. Wiss. R.  1983. 32 (1), 1983. pp. 41 - 46.\\
6.  Джрбашян М. М.,  Нерсесян А. Б. Дробные производные и задачи Коши для дифференциальных уравнений дробного порядка. Изв. АН АрмССР. Матем.,3:1 (1968), c. 3–28.\\
7.  Богатырева Ф.Т. Начальная задача для уравнения дробного порядка с постоянными коэффициентами. Вестник КРАУНЦ. Физ.-мат. науки, 2016, N. 5, c. 21–26\\
8. Псху А. В. Начальная задача для линейного обыкновенного дифференциального уравнения дробного порядка // Математический сборник. 2011. Т. 202. №4. c. 111-122\\
9.Karimov E., Ruzhansky M., Toshtemirov B.  Solvability of the boundary-value problem for a mixed equation involving hyper-Bessel fractional differential operator and bi-ordinal Hilfer fractional derivative. Mathematical Methods in the Applied Sciences. 41(1), 2023, pp. 54-77.\\
10. Restrepo, J. E., Suragan, D. (2021). Hilfer-type fractional differential equations with variable coefficients. Chaos, Solitons and Fractals, 150, 111146. doi:10.1016/j.chaos.2021.111146\\
11.Yuldashev T.K.,Kadirkulov B.J., Bandaliyev R.A. On a Mixed Problem for Hilfer Type Fractional Differential Equation with Degeneration.
Lobachevskii Journal of Mathematicsthis link is disabled, 2022, 43(1), pp. 263–274\\
12. B.Kh. Turmetov, B.J. Kadirkulov. On a problem for nonlocal mixed-type fractional order equation with degeneration.
Chaos, Solitons and Fractals,Volume 146, 2021, 110835, ISSN 0960-0779, https://doi.org/10.1016/j.chaos.2021.110835.\\
13. Smadiyeva A.G.   Well-posedness of the initial-boundary value problems for the
time-fractional degenerate diffusion equations. Bulletin of the Karaganda University. Mathematics series. 107(3), 2022, pp. 145-151.\\
14. Kilbas, Anatoly A.; Srivastava, Hari M.; Trujillo, Juan J. Theory and applications of fractional differential equations. North-Holland Mathematics Studies, 204. Elsevier Science B.V., Amsterdam, 2006.\\
15.Г.Бейтмен и А.Эрдейи. Высшие трансцендентные функции. Гипергеометрическая функция. Функция Лежандра. Издание второе. Изд. Наука, Москва, 1973, 296 С.

\end{document}